\documentclass[preprint,
aps,showpacs,preprintnumbers,amsmath,amssymb,floatfix,pre]{revtex4}

\usepackage{color,amsmath,graphicx,pst-text,%
amssymb,
            pst-grad,pstricks,pst-3d
}
\input epsf
\usepackage{epsfig}
\begin{document}

\newcommand{\GammaFn}[1]{\Gamma\left(#1\right)}
\newcommand{\Laplace}[2]{\mathcal{L}\left\lbrace #1\right\rbrace\left(#2\right)}
\newcommand{\Mellin}[2]{\mathcal{M}\left\lbrace #1\right\rbrace\left(#2\right)}
\newcommand{\Fourier}[2]{\mathcal{F}\left[ #1\right]\left(#2\right)}
\newcommand{\FourierC}[2]{\mathcal{F}_{c}\left[ #1\right]\left(#2\right)}
\newcommand{\Fourierb}[2]{\mathcal{F}\left\lbrace #1\right\rbrace\left(#2\right)}
\newcommand{\FourierCb}[2]{\mathcal{F}_{c}\left\lbrace #1\right\rbrace\left(#2\right)}
\newcommand{\Fracderiv}[1]{\mathcal{D}^{1-#1}}
\newcommand{\FracIntN}[1]{\frac{\partial^{#1}}{\partial t^{#1}}}
\newcommand{\FracIntM}[2]{\frac{\partial^{#1} #2 }{\partial t^{#1}}}
\newcommand{\FracderivN}[1]{\frac{\partial^{1-#1}}{\partial t^{1-#1}}}
\newcommand{\FracderivM}[2]{\frac{\partial^{1-#1} #2 }{\partial t^{1-#1}}}
\newcommand{\RLderivO}[1]{\frac{d^{#1}}{dt^{#1}}}
\newcommand{\RLderivOM}[2]{\frac{d^{#1} #2}{dt^{#1}}}

\title{Fractional Chemotaxis Diffusion Equations.}
\author{T.A.M. Langlands}
\email{t.langlands@usq.edu.au}
\affiliation{Department of Mathematics and Computing,
University of  Southern Queensland, Toowoomba Queensland 4350, Australia.}

\author{B.I. Henry}
\email{B.Henry@unsw.edu.au}
\affiliation{Department of Applied Mathematics, School of Mathematics,
University of New South Wales, Sydney NSW 2052, Australia.}

\date{\today}
\begin{abstract}
We introduce mesoscopic and macroscopic model equations of chemotaxis with anomalous subdiffusion  for modelling
chemically directed transport of biological organisms in changing chemical environments with diffusion hindered by
traps or macro-molecular crowding. 
 The mesoscopic models are formulated using Continuous Time Random Walk master equations and the macroscopic models are formulated with 
 fractional order differential equations. Different models are proposed depending on the timing of the chemotactic forcing. Generalizations
 of the models to  include linear reaction dynamics are also derived. Finally a Monte Carlo method for simulating anomalous subdiffusion with chemotaxis is introduced and simulation results are compared with numerical solutions of the  model equations. 
 The model equations developed here could be used to replace Keller-Segel type equations in biological systems with transport
 hindered by
traps, macro-molecular crowding or other obstacles.   
\end{abstract}
\maketitle
\section{Introduction}

Diffusion and
chemotaxis are fundamental to the motion of bacteria \cite{WA2004}, the directed motion of neutrophils in response to infection \cite{HD2004}, hypoxia stimulated angiogenesis \cite{OAMB2009}
and many other biological transport processes \cite{HD2004}. These transport processes can further be complicated by
 traps \cite{Saxton2007}, macromolecular crowding \cite{DV2008} or other obstacles resulting in
anomalous subdiffusion characterized by an ensemble averaged  mean square displacement of diffusing species, $\langle r^2(t)\rangle$, that scales
sublinearly in time, i.e., $\langle r^2(t)\rangle\sim t^\gamma$
with $0<\gamma<1$, 
\cite{Ghosh1994,Feder1996,Sa96,Sh97,S99,BWZW99,Simson1998,MK2000,Saxton2001,WHN2003,BF2005,OBSTVB2006}.
In this paper we introduce mesoscopic and macroscopic models for transport in biological systems with chemotaxis and anomalous subdiffusion.

The classic macroscopic model for the evolution of
a diffusing species, with concentration $n(x,t)$, in the presence of a  chemoattractant, with concentration $c(x,t)$, is the Keller-Segel model \cite{KS1971}
\begin{equation}
\frac{\partial n}{\partial t}=D\frac{\partial^2 n}{dx^2}-\chi\frac{\partial}{\partial x}\left(n\frac{\partial c}{\partial x}\right)\label{KS}
\end{equation}
where 
$D$ and $\chi$ denote the diffusion coefficient and the chemotactic coefficient respectively. In this model if the chemoattractant is removed
the evolution corresponds to standard Brownian diffusion with $\langle r^2(t)\rangle\sim t$.

Anomalous subdiffusion can be modelled as
fractional Brownian motion (fBm) \cite{A1976,Wang1990,Wang1994} or Continuous Time Random Walks (CTRWs) \cite{MW1965,SL1973} with long-tailed waiting-time densities
\cite{MK2000}.  Both of these models are non-Markovian and both exhibit the same sublinear
scaling for the ensemble averaged mean square displacement. However the second moment of the velocity scales differently in the two models
\cite{Lutz2001} and the time averaged mean square displacement differs from the ensemble averaged mean square displacements in the CTRW model,
but not in the fBm model \cite{HBMB2008}. Both possibilities should be considered when interpreting results from experiments using
 single particle tracking \cite{HBMB2008} or fluorescence recovery after photobleaching \cite{LK2008} and a simple 
 test has been devised for analysing experimental data to determine which model is most appropriate  \cite{MWBK2009}.

At the macroscopic level, anomalous subdiffusion can be modelled through a modified diffusion equation
\begin{equation}
\frac{\partial C}{\partial t}={\mathcal D}(\gamma, t)\nabla^2 C
\end{equation}
with the diffusion constant replaced by a fractional temporal operator. In the case of
 fractional Brownian motion (fBm) this operator is given by
  \cite{Wang1990,Wang1994}
\begin{equation}
{\mathcal D}_I(\gamma, t)=D(\gamma)\gamma t^{\gamma-1}\label{DfBm}
\end{equation}
In the Continuous Time Random Walk (CTRW) model \cite{MW1965}, with power law waiting times
 \cite{MK2000}, the fractional temporal operator is given by
\begin{equation}
{\mathcal D}_{II}(\gamma,t)=D(\gamma)\frac{\partial^{1-\gamma}}{\partial t^{1-\gamma}}\label{DCTRW}
\end{equation}
where
$D(\gamma)$ is a generalized diffusion coefficient
 with units of
 $m^2s^{-\gamma}$ and
 \begin{equation}
 \frac{\partial^{1-\gamma}}{\partial t^{1-\gamma}}Y(t)
 =\frac{1}{\Gamma(\gamma)}\frac{\partial}{\partial t}\int_0^t
 \frac{Y(t')}{(t-t')^{1-\gamma}}dt'
 \end{equation}
 defines the Riemann-Liouville fractional derivative
 of order $1-\gamma$. 
There have been various attempts to modify the fractional macroscopic equations to include force fields and reactions
\cite{MK2000,HLS2010}. 

The Fokker-Planck equation for diffusion in a force field can readily be generalized by replacing the diffusion coefficient with a time dependent
fractional operator as above. This has been justified within the framework of
 CTRWs, for force fields that
vary in space but not time \cite{Barkai2000,MK2000} and for force fields that vary in time but not
space \cite{SK2006}.  However these derivations do not extend to the more general  case of anomalous subdiffusion in a general external force field $f(x,t)$ that varies in both  time and space. Two obvious possible generalizations in this case are \cite{Sokolov2001}
\begin{equation}
 \label{fokkerplanck}
 \frac{\partial n}{\partial t} = \frac{\partial^{1-\gamma}}{\partial t^{1-\gamma}}
 D_\gamma \nabla^2 n -\frac{1}{\eta_\gamma}  \frac{\partial^{1-\gamma}}{\partial t^{1-\gamma}}
 \nabla\left(f(x,t) n(x,t)\right)
\end{equation}
and \cite{HPGH2007,WMW2008,HPGH2009}
\begin{equation}
 \label{dispersal}
 \frac{\partial n}{\partial t} = \frac{\partial^{1-\gamma}}{\partial t^{1-\gamma}}
 D_\gamma \nabla^2 n -\frac{1}{\eta_\gamma}  \nabla\left(
f(x,t)  \frac{\partial^{1-\gamma}}{\partial t^{1-\gamma}}
  n(x,t)\right).
\end{equation}
If the force field is purely space dependent then the two models are equivalent and the solution is 
 time subordinated to the concentration
of diffusing species in the standard Fokker-Planck equation. This temporal  subordination is not physically appropriate for time dependent 
external force fields \cite{HPGH2007}. However an alternate formulation using an Ito stochastic differential equation has been proposed with a modified subordination in which the force varies in real time rather than the random time \cite{WMW2008}.   In chemotaxis  there may be a physical link between the time
scale of the diffusion and the time scale of the effective force field since the latter
depends on the concentration of another
diffusing species. Similarly in the fractional Nernst-Planck equation considered in \cite{HLW2008,LHW2009}
the force field from the membrane potential
 depends on concentrations of the diffusing species.

In Section II we introduce four different models of chemotaxis with anomalous subdiffusion. The different models are characterized
by differences in the nature 
of the anomalous diffusion (fBm or power law CTRWs), and differences in the details
of the underlying random walk processes.  In Section III  numerical solutions of the associated discrete space equations are obtained for each model. The numerical results are compared with Monte Carlo random walk simulations, with chemotactic forcing, on the same grid and using the same parameters. 
Differences between the model results are discussed in Section IV.

\section{Fractional Chemotaxis Diffusion Models}
\label{Models}
\subsection{Model I}

To model chemotaxis with fractional Brownian motion we consider
 an {\em ad hoc} model in which we replace both the diffusion coefficient and the chemotactic coefficient by fractional temporal operators as in  Eq. (\ref{DfBm}). This yields  
\begin{equation}
\label{ModelI}
\frac{\partial n}{\partial t} = \gamma  t^{\gamma-1}\left[D_{\gamma}\frac{\partial^2 n}{\partial x^2} -\chi_\gamma\frac{\partial}{\partial x}\left(n\frac{\partial c}{\partial x}  \right) \right] 
\end{equation} 
where $\gamma$ is the anomalous diffusion exponent, $D_{\gamma}$ is the anomalous diffusion coefficient (with units $m^2s^{-\gamma}$),
and $\chi_{\gamma}$ is the analogous anomalous chemotaxis coefficient. 
This model equation reduces to the standard Keller-Segel chemotaxis equation, Eq. (\ref{KS}), when $\gamma=1$.

\subsection{Model II}
A simple model for chemotaxis with fractional diffusion from power law CTRWs starts with the master equation
 \begin{equation}
\label{ModelII:1}
n_{i}(t) = n_i(0) \Phi(t) + \int\limits_0^{t} \left\lbrace 
p_r(x_{i-1},t^{'}) n_{i-1}(t^{'})+
p_l(x_{i+1},t^{'}) n_{i+1}(t^{'}) \right\rbrace \psi(t-t^{'}) \:dt^{'} 
\end{equation} 
where $\psi(t)$ is a (power law) waiting time density,
\begin{equation}
\label{ModelII:surv}
\Phi(t) = \int\limits_{t}^{\infty} \psi(t^{'}) \; dt^{'}
\end{equation}
 is the corresponding survival probability,
and $p_r(x,t)$ and $p_l(x,t)$ are the probabilities of jumping  from $x$ to the adjacent grid point to the right and left directions respectively.
 These probabilities are dependent on the
chemoattractant concentrations, $c(x,t)$, at the neighbouring points of the point $x$ at time $t$. The master equation, Eq. (\ref{ModelII:1}), is a continuous time representation of the transition probability law in \cite{Stevens2000}.

Following Stevens \cite{Stevens2000}, the probabilities of jumping to the left or right direction are based on the proportion of the chemoattractant on either side of the current point via
\begin{equation}
\label{ModelII:pleft}
p_l(x_{i},t) = \frac{v(x_{i-1},t)}{v(x_{i-1},t)+v(x_{i+1},t)},
\end{equation} 
and
\begin{equation}
\label{ModelII:pright}
p_r(x_{i},t) = \frac{v(x_{i+1},t)}{v(x_{i-1},t)+v(x_{i+1},t)},
\end{equation} 
where $v(x,t)$ is a sensitivity function that depends on the concentration of the chemoattractant:
\begin{equation}
\label{ModelII:v}
v(x,t) = \exp{\left(\beta\,c(x,t)\right)}.
\end{equation} 
Note that with the above we have
\begin{equation}
p_l(x_i,t)+p_r(x_i,t)=1\label{normp}
\end{equation}
and
\begin{equation}
p_l(x_i,t)-p_r(x_i,t)=\frac{e^{\beta c(x_{i-1},t)}-e^{\beta c(x_{i+1},t)}}
{e^{\beta c(x_{i-1},t)}+e^{\beta c(x_{i+1},t)}}.\label{diffp}
\end{equation}

Using the notation $\Laplace{f(t)}{s}$ or $\widehat f(s)$ to denote the Laplace transform with respect to time of a function $f(t)$ we have
the Laplace transform of Eq. (\ref{ModelII:1}), 
 \begin{equation}
\label{ModelII:2}
\widehat{n}_{i}(s) = n_{i}(0) \widehat{\Phi}(s) +  
\left\lbrace \Laplace{p_r(x_{i-1},t) n_{i-1}(t)}{s}
+
\Laplace{p_l(x_{i+1},t) n_{i+1}(t)}{s}\right\rbrace  \widehat{\psi}(s).
\end{equation} 
Using the identity
\begin{equation}
\label{Lap:surv}\widehat{\Phi}(s) =(1-\widehat{\psi}(s))/s,
\end{equation}
which follows from the Laplace transform of (\ref{ModelII:surv}),
we have
 \begin{multline}
\label{ModelII:3}
s\widehat{n}_{i}(s) - n_{i}(0)=\frac{\widehat{\psi}(s)}{\widehat{\Phi}(s)}  
\left\lbrace -\widehat{n}_{i}(s)
+\Laplace{p_r(x_{i-1},t) n_{i-1}(t)}{s}
+\Laplace{p_l(x_{i+1},t) n_{i+1}(t)}{s}\right\rbrace. 
\end{multline} 

We now consider a heavy-tailed waiting-time density which behaves for long-times as
\begin{equation}
\label{heavy:tail} \psi(t)\sim \frac{\kappa}{\tau}\left(\frac{t}{\tau}\right)^{-1-\gamma} 
\end{equation}
where $\gamma$  is
the anomalous exponent, $\tau$ is the characteristic waiting-time, and $\kappa$ is a dimensionless constant.
Using a Tauberian (Abelian) theorem \cite{Margolin2004} we can write the Laplace transform for this density function as (for small $s$)
\begin{equation}
\label{heavy:tail:laplace} \widehat{\psi}(s)\sim 1 -  \frac{\kappa\GammaFn{1-\gamma}}{\gamma}\left(s\tau\right)^{\gamma}. 
\end{equation}
Using Eq.~(\ref{Lap:surv}), we then find the corresponding asymptotic form for the survival probability
\begin{equation}
\label{surv:asymp:laplace} \widehat{\Phi}(s)\sim   \frac{\kappa\GammaFn{1-\gamma}}{\gamma}\tau^{\gamma} s^{\gamma-1}  
\end{equation}
and the ratio  
\begin{equation}
\label{L:Ratio:ht}
 \frac{\widehat{\psi}(s)}{\widehat{\Phi}(s)} \sim A_\gamma
\frac{s^{1-\gamma}}{\tau^{\gamma}}
\end{equation}
where $$A_\gamma=\frac{\gamma}{\kappa\GammaFn{1-\gamma}}.$$

Specific cases of waiting time densities are the Mittag-Leffler density \cite{SGMR03}
\begin{equation}
\label{ModelII:Mittag:psi}
\psi(t) = -\frac{d}{dt} E_{\gamma}\left(-\left(\frac{t}{\tau}\right)^\gamma\right),
\end{equation} 
where $E_{\gamma}(z)$ is the Mittag-Leffler function \cite{Podlubny1999},
and
the Pareto  law used by \cite{Yuste2004} 
\begin{equation}
\label{Pareto:psi}
\psi\left(t\right) = \frac{{\gamma}/{\tau}}{\left(1+t/\tau\right)^{1+\gamma}}.
\end{equation} 
The corresponding values for $A_\gamma$ can be shown to be
\begin{equation}
A_{\gamma}= 1
\quad \mathrm{and} \quad
\label{Agamma:Pd}
A_{\gamma}= \frac{1}{\GammaFn{1-\gamma}}
\end{equation}
for (\ref{ModelII:Mittag:psi})
and (\ref{Pareto:psi}) respectively.
Note the ratio in Eq.(\ref{L:Ratio:ht}) is only valid  long-times for the Pareto density, Eq.~(\ref{Pareto:psi}), whilst it is exact for the Mittag-Leffler density for all times. In addition, if $\gamma=1$ we do not use  Eq.~(\ref{Pareto:psi}) but instead use Eq.~(\ref{ModelII:Mittag:psi}).

With Eq. (\ref{L:Ratio:ht}), Eq. (\ref{ModelII:3}) now
becomes 
 \begin{multline}
\label{ModelII:4}
s\widehat{n}_{i}(s) - n_{i}(0)=  
\frac{A_{\gamma} s^{1-\gamma}}{\tau^\gamma}\left\{-\widehat{n}_{i}(s)+\Laplace{p_r(x_{i-1},t) n_{i-1}(t)}{s}+\Laplace{p_l(x_{i+1},t) n_{i+1}(t)}{s}\right\}.
\end{multline} 
Noting that the Laplace Transform of a Riemann-Liouville fractional derivative of order $\alpha$, where $0<\alpha\le 1$, is given by \cite{Podlubny1999}
\begin{equation}
\label{Laplace:FracDeriv} 
\Laplace{\RLderivOM{\alpha}{f(t)}}{s} = s^{\alpha} \widehat{f}(s)- 
\left[\RLderivOM{\alpha-1}{f(t)}\right|_{t=0}
\end{equation} 
we can invert the Laplace transforms in
Eq.(\ref{ModelII:4}) to obtain
\begin{equation}
\label{ModelII:8}
\frac{d n_i}{dt} = \frac{A_{\gamma}}{\tau^\gamma} \RLderivO{1-\gamma} \left\lbrace  -n_i(t)+
p_r(x_{i-1},t) n_{i-1}(t) +p_l(x_{i+1},t) n_{i+1}(t)\right\rbrace
\end{equation} 
where we have ignored the last term in Eq. (\ref{Laplace:FracDeriv}). Numerical solutions of this discrete space fractional differential equation for Model II are considered in Section III.  

The spatial continuum limit of Model II can be obtained in the usual way 
by setting  $x_{i}=x$ and $x_{i\pm 1}=x\pm \Delta x$ and carrying out Taylor series expansions in $x$.
Retaining terms to order $(\Delta x)^2$ and using the normalization
$p_l(x,t)+p_r(x,t)=1$ we first find that
 \begin{eqnarray}
 -n_i(t)+
p_r(x_{i-1},t) n_{i-1}(t) +p_l(x_{i+1},t) n_{i+1}(t)&=&
\Delta x\frac{\partial}{\partial x}\left(n(x,t)[p_l(x,t)-p_r(x,t)]\right)\nonumber\\
& &+\frac{\Delta x^2}{2}\frac{\partial^2}{\partial x^2}n(x,t).\label{result1}
\end{eqnarray}
This simplifies further after carrying out Taylor series expansions in Eq.(\ref{diffp}), to arrive at
\begin{eqnarray}
p_l(x,t)-p_r(x,t)&\approx&\frac{e^{-\beta\Delta x\frac{\partial c}{\partial x}}
-e^{\beta\Delta x\frac{\partial c}{\partial x}}}{e^{-\beta\Delta x\frac{\partial c}{\partial x}}
+e^{\beta\Delta x\frac{\partial c}{\partial x}}},\nonumber\\
&=&-\tanh\left(\beta\Delta x\frac{\partial c}{\partial x}\right),\nonumber\\
&\approx&-\beta\Delta x\frac{\partial c}{\partial x}.\label{result2}
\end{eqnarray}

We can now combine the results in
Eq.(\ref{result1}) and Eq.(\ref{result2}) with Eq. (\ref{ModelII:8}) to obtain
\begin{equation}
\label{ModelII:eqn}
\frac{\partial n}{\partial t} \simeq \FracderivN{\gamma} \left[ \frac{A_{\gamma}\Delta x^2}{2\tau^\gamma}\frac{\partial^{2}n}{\partial x^2}(x,t)-\frac{A_{\gamma}\beta\Delta x^2}{\tau^\gamma}\frac{\partial}{\partial x}\left(\frac{\partial c(x,t)}{\partial x} n(x,t)\right) \right]+O(\Delta x^4)
\end{equation} 
and then taking the limit $\Delta x\rightarrow 0$  and $\tau\rightarrow 0$, with
\begin{equation}
\label{Diffusion:coeff}
D_{\gamma} = \frac{A_{\gamma}\Delta x^2}{2\tau^\gamma}
\end{equation} 
and 
\begin{equation}
\label{Chemo:coeff}
\chi_{\gamma} = \frac{A_{\gamma}\beta\Delta x^2}{\tau^\gamma},
\end{equation} 
we have 
\begin{equation}
\label{ModelII:eqn:final}
\frac{\partial n}{\partial t} = \FracderivN{\gamma} \left[ 
D_{\gamma}\frac{\partial^{2}n(x,t)}{\partial x^2}-\chi_{\gamma}\frac{\partial}{\partial x}\left(\frac{\partial c(x,t)}{\partial x} n(x,t)\right) \right].
\end{equation} 
Equation (\ref{ModelII:eqn:final}) provides a useful approximation
for the space and time evolution of the concentration
of an anomalously diffusing species that is chemotactically attracted by another species. In the CTRW  master equation, Eq.(\ref{ModelII:1}) for this model the probabilities to jump left or right are determined at the start of the waiting times. We will consider
another model in Section D where the probabilities to jump left or right are determined at the end of the waiting times. But first, in the next section, we consider an alternate formulation, based on a generalized master equation approach.

\subsection{Model III}
In this section we follow the  generalised  master equation approach of \cite{Chechkin2005,SK2007,SK2006}, extended to take into account the effect of the chemoattractant.
To begin we write the balance equation for  the concentration of particles, $n$, at the site $i$
\begin{equation}
\label{Model:III:Balance}
\frac{dn_i(t)}{dt} = J_{i}^{+}(t)-J_{i}^{-}(t),
\end{equation} 
where $J_{i}^{\pm}$ are the gain (+) and loss (-) fluxes at the site $i$.
We also  have the conservation equation for the arriving flux of particles at the point $i$ 
given by the flux of particles either leaving the site $i-1$  and jumping to the right or leaving the site $i+1$ that move to the left:
\begin{equation}
\label{Model:III:Gain:flux}
J_{i}^{+}(t)=  p_{r}(x_{i-1},t)J_{i-1}^{-}(t)+ p_{l}(x_{i+1},t)J_{i+1}^{-}(t)
\end{equation} 
where $p_{l}(x,t)$ and $p_{r}(x,t)$ are given in Eqs.~(\ref{ModelII:pleft}) and (\ref{ModelII:pright}).
We can combine Eq.~(\ref{Model:III:Gain:flux})  and Eq.~ (\ref{Model:III:Balance}) to obtain an evolution law for the concentration purely
in terms of the loss flux, viz;
\begin{equation}
\label{Model:III:Balance:2}
\frac{dn_i(t)}{dt} = p_{r}(x_{i-1},t)J_{i-1}^{-}(t)+ p_{l}(x_{i+1},t)J_{i+1}^{-}(t)-J_{i}^{-}(t).
\end{equation} 
The loss flux at the site $i$  is given by
\begin{equation}
\label{Model:III:Loss:flux}
J_{i}^{-}(t)= \psi(t)n_{i}(0)+\int\limits_{0}^{t} \psi(t-t^{'}) J_{i}^{+}(t^{'}) \; dt^{'}.
\end{equation} 
The first  term represents those particles that were either originally at $i$ at $t=0$ and wait until time $t$ when they leave. The second term
represents particles that arrived at some earlier time $t^{'}$ and wait until time $t$ to leave. Here $\psi(t)$ is the 
usual waiting-time density used in Model II.
We
can combine Eq.~(\ref{Model:III:Balance}) and Eq.~(\ref{Model:III:Loss:flux}) to obtain
\begin{equation}
\label{Model:III:Loss:flux:sub}
J_{i}^{-}(t)= \psi(t)n_{i}(0)+\int\limits_{0}^{t} \psi(t-t^{'}) \left[J_{i}^{-}(t^{'})+\frac{dn_i(t^{'})}{dt}\right] \; dt^{'}
\end{equation} 
and then we can solve for the loss flux using Laplace transform methods.
The Laplace transform of Eq.(\ref{Model:III:Loss:flux:sub}) with respect to time yields
\begin{equation}
\label{Model:III:Loss:flux:lapl}
\widehat{J}_{i}^{-}(s)= \widehat{\psi}(s)n_{i}(0)+ \widehat{\psi}(s) \left[\widehat{J}_{i}^{-}(s)+s\widehat{n}_i(s)-n_{i}(0)\right] 
\end{equation} 
which simplifies further as
\begin{equation}
\label{Model:III:Loss:flux:lapl:2}
\widehat{J}_{i}^{-}(s)= \frac{\widehat{\psi}(s)}{\widehat{\phi}(s)}\widehat{n}_i(s).
\end{equation} 
Now using the approximation in Eq.~(\ref{L:Ratio:ht}) for a heavy-tailed waiting-time density 
 and inverting the Laplace transform we have
\begin{equation}
\label{Model:III:Loss:flux:invert}
{J}_{i}^{-}(t)= \frac{A_{\gamma}}{\tau^\gamma}\RLderivOM{1-\gamma}{n_i(t)}.
\end{equation} 

We now substitute the expression for the loss flux, Eq.(\ref{Model:III:Loss:flux:invert})
back into the balance equation, Eq.(\ref{Model:III:Balance:2}),
to obtain
\begin{equation}
\label{Model:III:Balance:3}
\frac{dn_i(t)}{dt} =\frac{A_{\gamma}}{\tau^\gamma}\left\lbrace
                    p_{r}(x_{i-1},t)\RLderivOM{1-\gamma}{n_{i-1}(t)}+
                    p_{l}(x_{i+1},t)\RLderivOM{1-\gamma}{n_{i+1}(t)}-
                    \RLderivOM{1-\gamma}{n_i(t)}\right\rbrace.
\end{equation} 
Numerical solutions of this discrete space fractional differential equation for Model III are considered in Section III.

The continuous space representation
of Eq.(\ref{Model:III:Balance:3}) is found by setting
 $x_{i}=x$ and $x_{i\pm 1}=x\pm \Delta x$ so that
\begin{multline}
\label{Model:III:Balance:4}
\frac{\partial n(x,t)}{\partial t} =\frac{A_{\gamma}}{\tau^\gamma}\left\lbrace
                          p_{r}(x-\Delta x,t)\FracderivM{\gamma}{n(x-\Delta x,t)}\right.\\
                          \left.
                          +p_{l}(x+\Delta x,t)\FracderivM{\gamma}{n(x+\Delta x ,t)}-
                           \FracderivM{\gamma}{n(x,t)}\right\rbrace.
\end{multline} 
The continuum limit representation can then be found by carrying out Taylor series expansions about $x$,  similar to the steps used to reduce  (\ref{ModelII:8}) to (\ref{ModelII:eqn:final}). This results in the equation
\begin{equation}
\label{ModelIII:eqn:final}
\frac{\partial n}{\partial t} = \FracderivN{\gamma}  
D_{\gamma}\frac{\partial^{2}n(x,t)}{\partial x^2}-
\chi_{\gamma}\frac{\partial}{\partial x}\left(\frac{\partial c(x,t)}{\partial x} \FracderivM{\gamma}{n(x,t)}\right).
\end{equation} 

Model II and Model III are  similar to the fractional Fokker Planck  equations, Eq.(\ref{fokkerplanck}) and Eq.(\ref{dispersal}) respectively with forcing from the chemotactic gradient $\frac{\partial c(x,t)}{\partial x}$.

\subsection{Model IV}
We now  re-consider the master equation for the CTRW model but with the jump probabilities calculated after the particle has waited and immediately prior to jumping. The master equation in this case is given by
 \begin{equation}
\label{ModelIV:1}
n_{i}(t) = n_{i}(0) \Phi(t) + p_r(x_{i-1},t) \int\limits_0^{t} 
n_{i-1}(t^{'}) \psi(t-t^{'})\:dt^{'}+
p_l(x_{i+1},t)\int\limits_0^{t} n_{i+1}(t^{'})  \psi(t-t^{'}) \:dt^{'}.
\end{equation} 
It is convenient to introduce the auxiliary function
\begin{equation}
m_i(t)=\int_0^t n_i(t^{'})\psi(t-t^{'})\:dt^{'},\label{mterm}
\end{equation}
which has the Laplace transform
\begin{equation}
\widehat{m}_i(s)=\widehat{n}_i(s)\widehat\psi(s).
\end{equation}
The Laplace transform of Eq.(\ref{ModelIV:1}) with respect to time
can then be written as
 \begin{equation}
\label{ModelIVL}
\left(s\widehat{n}_{i}(s) - n_{i}(0)\right)\widehat\Phi(s)= - \widehat\psi(s)
\widehat{n}_{i}(s)+\Laplace{p_r(x_{i-1},t) m_{i-1}(t)}{s}+\Laplace{p_l(x_{i+1},t) m_{i+1}(t)}{s},
\end{equation} 
and after the inverse  Laplace transform, 
\begin{equation}
\int_0^t \frac{\partial n_i}{\partial t^{'}}\Phi(t-t^{'})\, dt^{'}
=m_i(t)+p_r(x_{i-1},t) m_{i-1}(t)+p_l(x_{i+1},t) m_{i+1}(t).
\end{equation}
Proceeding to the continuum limit with Taylor series expansions about
$x$, similar to the steps in Model II and Model III, we obtain
\begin{equation}
\int\limits_0^{t} \Phi(t-t^{'}) \frac{\partial n(x,t^{'})}{\partial t} \:dt^{'}  \simeq  \frac{\Delta x^2}{2}
\frac{\partial^2 m(x,t)}{\partial x^2}   
-\Delta{x}^2\beta\frac{\partial}{\partial x}\left(\frac{\partial c(x,t)}{\partial x} 
m(x,t) \right)+O(\Delta x^4),
\end{equation} 
and then using the auxiliary function definition in Eq.(\ref{mterm}) we find
\begin{multline}
\label{ModelII:4:alt}
\int\limits_0^{t} \Phi(t-t^{'}) \frac{\partial n(x,t^{'})}{\partial t} \:dt^{'}  \simeq  \frac{\Delta x^2}{2}
\int\limits_0^{t} \frac{\partial^2 n(x,t^{'})}{\partial x^2}   \psi(t-t^{'}) \:dt^{'}\\
-\Delta{x}^2\beta\frac{\partial}{\partial x}\left(\frac{\partial c(x,t)}{\partial x}\int\limits_0^{t} 
n(x,t^{'}) \psi(t-t^{'})\:dt^{'}\right)+O(\Delta x^4).
\end{multline} 

Asymptotic expressions for the convolution integrals in Eq.(\ref{ModelII:4:alt}) can be obtained by considering asymptotic expansions in Laplace space and then inverting.
Thus we now  consider the terms
\begin{eqnarray}
\label{ModelII:alt:asym:t1}
\Laplace{\int\limits_0^{t} \Phi(t-t^{'}) \frac{\partial n(x,t^{'})}{\partial t} \:dt^{'}}{s}  &=&\widehat{\Phi}(s) \Laplace{\frac{\partial n(x,t)}{\partial t}}{s}, \\
\label{ModelII:alt:asym:t2}
\Laplace{\int\limits_0^{t} \frac{\partial^2 n(x,t^{'})}{\partial x^2}   \psi(t-t^{'}) \:dt^{'}}{s} &=&
 \widehat{\psi}(s) \frac{\partial^2 \widehat{n}(x,s)}{\partial x^2},\\
\label{ModelII:alt:asym:t3}
\Laplace{\int\limits_0^{t} n(x,t^{'}) \psi(t-t^{'})\:dt^{'}}{s} &=& \widehat{\psi}(s) \widehat{n}(x,s).
\end{eqnarray}
For long times (small $s$) we have, with the use of Eqs.~(\ref{heavy:tail:laplace}) and (\ref{surv:asymp:laplace})
\begin{eqnarray}
\label{ModelII:alt:asym:Phi:lg}
\widehat{\Phi}(s)
 &\simeq& \frac{s^{\gamma-1}\tau^\gamma}{A_{\gamma}} +O(s^{2\gamma-1} ), \\
\label{ModelII:alt:asym:psi:lg}
\widehat{\psi}(s) &\simeq & 1-\frac{\left(s\tau\right)^{\gamma}}{A_{\gamma}}+O(s^\gamma).
\end{eqnarray}
Using these expansions in  Eqs.~(\ref{ModelII:alt:asym:t1})~--~(\ref{ModelII:alt:asym:t3}) and taking the inverse Laplace transforms to replace the convolution integrals in  Eq.~(\ref{ModelII:4:alt}) we obtain 
\begin{equation}
\label{ModelII:4:alt:lg}
\frac{\tau^\gamma}{A_\gamma} \FracIntN{\gamma-1}\frac{\partial n(x,t)}{\partial t} \simeq \frac{\Delta x^2}{2}
\frac{\partial^2 n(x,t)}{\partial x^2} 
-\Delta{x}^2\beta\frac{\partial}{\partial x}\left(\frac{\partial c(x,t)}{\partial x} 
n(x,t)\right)+O(\Delta x^4),
\end{equation} 
and in the limit $\Delta x\rightarrow 0$ and $\tau \rightarrow 0$
\begin{equation}
\label{ModelII:4:alt:lg:2}
\frac{\partial n(x,t)}{\partial t} = \FracIntN{1-\gamma} \left[D_\gamma
\frac{\partial^2 n(x,t)}{\partial x^2} 
-\chi_\gamma\frac{\partial}{\partial x}\left(\frac{\partial c(x,t)}{\partial x} 
n(x,t)\right)\right]
\end{equation} 
as previously in Eq.~(\ref{ModelII:eqn:final}).

Conversely, for short times (large $s$) we have
\begin{eqnarray}
\label{ModelII:alt:asym:Phi:sh}
\widehat{\Phi}(s)
 &\simeq& \frac{1}{s} +O(s^{-\nu_\gamma-1}), \\
\label{ModelII:alt:asym:psi:sh}
\widehat{\psi}(s) &\simeq & \frac{B_{\gamma}s^{-\nu_\gamma}}{\tau^{\nu_\gamma}}+O(s^{-2\nu_\gamma}),
\end{eqnarray}
where $\nu_\gamma =\gamma$ and $B_{\gamma}=1$ if use the Mittag-Leffler (\ref{ModelII:Mittag:psi}) and $\nu_\gamma =1$ and $B_{\gamma}=\gamma$
if we use the Pareto (\ref{Pareto:psi}) density.
The resulting equation for Eq.~(\ref{ModelII:4:alt}) becomes for short-times 
\begin{equation}
\label{ModelII:4:alt:sh}
\frac{\partial n(x,t)}{\partial t} =  D^{*}_\gamma\FracIntN{1-\nu_\gamma} \frac{\partial^2 n(x,t)}{\partial x^2} 
-\chi^{*}_\gamma\frac{\partial}{\partial x}\frac{\partial}{\partial t}\left(\frac{\partial c(x,t)}{\partial x} 
\FracIntM{-\nu_\gamma}{n(x,t)}\right)
\end{equation} 
with the modified coefficients
\begin{equation}
\label{Diffusion:coeff:mod}
D_{\gamma}^{*} = \frac{B_{\gamma}\Delta x^2}{2\tau^{\nu_\gamma}}
\end{equation} 
and 
\begin{equation}
\label{Chemo:coeff:mod}
\chi_{\gamma}^{*} = \frac{B_{\gamma}\beta\Delta x^2}{\tau^{\nu_\gamma}}.
\end{equation} 
In the case of the Mittag-Leffler density (\ref{ModelII:Mittag:psi}), 
the short-time equation can be simplified to
\begin{multline}
\label{ModelII:4:alt:sh:ML}
\frac{\partial n(x,t)}{\partial t} =  D_\gamma^{*}\FracIntN{1-\gamma} \frac{\partial^2 n(x,t)}{\partial x^2} 
-\chi_\gamma\frac{\partial}{\partial x}\left(\frac{\partial c(x,t)}{\partial x}\FracIntM{1-\gamma}{n(x,t)}\right)
\\ -\chi_\gamma^{*}\frac{\partial}{\partial x}\left(\frac{\partial^2 c(x,t)}{\partial x\partial t} 
\FracIntM{-\gamma}{n(x,t)}\right).
\end{multline}
Equation (\ref{ModelII:4:alt:sh:ML}) is similar to Model III's  governing equation given  by Eq.~(\ref{ModelIII:eqn:final}) if we note that the last term which will be close to zero near $t=0$. Conversely if we use the Pareto density (\ref{Pareto:psi}) then   Eq.~(\ref{ModelII:4:alt:sh}) becomes
\begin{multline}
\label{ModelII:4:alt:sh:PD}
\frac{\partial n(x,t)}{\partial t} =  \gamma D_1 \frac{\partial^2 n(x,t)}{\partial x^2} 
-\gamma \chi_1\frac{\partial}{\partial x}\left(\frac{\partial c(x,t)}{\partial x} n(x,t)\right)
\\ -\gamma\chi_1\frac{\partial}{\partial x}\left(\frac{\partial^2 c(x,t)}{\partial x\partial t} 
\FracIntM{-1}{n(x,t)}\right)
\end{multline}
which is similar to the standard chemotaxis equation (\ref{KS}) (if we again consider the last term to be small) except for the presence of the $\gamma$ term. This term will, in effect, slow the 
initial temporal behaviour of the solution (linear rescaling). 

Note that if we use the Mittag-Leffler density in (\ref{ModelII:Mittag:psi}) then we see that the mesoscopic equation
Eq.~(\ref{ModelIV:1}) bridges the gap between Model II and Model III. At short times it recovers Model III, Eq.(\ref{ModelIII:eqn:final}),  whereas at long times it recovers Model II (\ref{ModelII:eqn:final}). The  latter shows that for long times, compared to the characteristic time $\tau$, there is no difference between evaluating the
chemotactic probabilities at the time before or after the particle waits.

\section{Fractional Chemotaxis Reaction-Diffusion Models}
In this section we consider extensions of the CTRW based fractional chemotaxis diffusion models to incorporate reactions. In the absence of chemotaxis, extensions of CTRW based fractional diffusion models to include linear reactions  were derived in \cite{HLW2006,Sokolov2006} and  extensions to include nonlinear reactions were derived in \cite{Yadav2006,Fedotov2010}.

\subsection{Model II}
Following the approach in \cite{HLW2006} we can incorporate reactions in CTRW models by increasing or decreasing the concentration of particles
during the waiting times by an amount proportional to the evolution operator for the reaction dynamics.  The master equation for Model II with linear reaction dynamics incorporated in this way becomes
 \begin{equation}
\label{ModelII:1:lr}
n_{i}(t) =e^{kt} n_{i}(0) \Phi(t) + \int\limits_0^{t} \left\lbrace 
p_r(x_{i-1},t^{'}) n_{i-1}(t^{'})
+
p_l(x_{i+1},t^{'}) n_{i+1}(t^{'}) \right\rbrace e^{k\left(t-t^{'}\right)}\psi(t-t^{'}) \:dt^{'} 
\end{equation} 
where $k$ is the per capita rate gain ($k>0$) or loss ($k<0$) of particles.
Again Laplace transform methods can be used to convert the integral equation representation into a (fractional) differential equation. The Laplace transform of Eq.~(\ref{ModelII:1:lr}) with respect to time yields
 \begin{multline}
\label{ModelII:2:lr}
\widehat{n}_{i}(s) = n_{i}(0) \widehat{\Phi}(s-k) +  
\Laplace{p_r(x_{i-1},t) n_{i-1}(t)}{s}\widehat{\psi}(s-k)
\\+
\Laplace{p_l(x_{i+1},t) n_{i+1}(t)}{s} \widehat{\psi}(s-k),
\end{multline} 
and after rearranging we find
 \begin{multline}
\label{ModelII:3:lr}
s\widehat{n}_{i}(s) - n_{i}(0)=  k\widehat{n}_{i}(s)
+\frac{\widehat{\psi}(s-k)}{\widehat{\Phi}(s-k)}
\left\lbrace
-\widehat{n}_{i}(s)+
\Laplace{p_r(x_{i-1},t) n_{i-1}(t)}{s}
\right.\\\left.+
\Laplace{p_l(x_{i+1},t) n_{i+1}(t)}{s}\right\rbrace,
\end{multline} 
where we have used (\ref{Lap:surv}).
With the result in Eq.(\ref{L:Ratio:ht}) we can inverting the Laplace transform to obtain 
\begin{multline}
\label{ModelII:8:lr}
\frac{dn_{i}}{dt} =  e^{kt}\frac{A_{\gamma}}{\tau^\gamma}\RLderivO{1-\gamma} \left(e^{-kt}\left[  -n_{i}(t)+
p_r(x_{i-1},t) n_{i-1}(t)+ p_l(x_{i+1},t) n_{i+1}(t)\right]\right)
+k n_{i}(t)
\end{multline} 
where the Riemann-Liouville fractional derivative has been replaced by a modified fractional derivative \cite{HLW2006,Sokolov2006}.

The  continuum limit, found by taking Taylor series expansions about $x$,  is
\begin{equation}
\label{ModelII:eqn:final:lr}
\frac{\partial n}{\partial t} = e^{kt} \FracderivN{\gamma} \left( e^{-kt}\left[ 
D_{\gamma}\frac{\partial^{2}n(x,t)}{\partial x^2}-\chi_{\gamma}\frac{\partial}{\partial x}\left(\frac{\partial c(x,t)}{\partial x} n(x,t)\right) \right]\right)+k n(x,t).
\end{equation} 
We note if $n(x,t)$ is not self-chemotactic then the solution of  (\ref{ModelII:eqn:final:lr}) is given by
$n(x,t)=e^{kt}y(x,t)$   where $y(x,t)$ is the solution of (\ref{ModelII:eqn:final}) with $y$ replacing $n$. 

\subsection{Model III}
To incorporate reactions in  Model III  we start by modifying Eq.~(\ref{Model:III:Balance}) to
\begin{equation}
\label{Model:III:Balance:lr}
\frac{dn_i(t)}{dt} = J_{i}^{+}(t)-J_{i}^{-}(t)+ k n_i(t)
\end{equation} 
where $k$, again, is the per capita rate gain ($k>0$) or loss ($k<0$) of particles.
We also  modify the expression for the loss flux, $ J_{i}^{-}(t)$,  in (\ref{Model:III:Loss:flux}) to
\begin{equation}
\label{Model:III:Loss:flux:lr}
J_{i}^{-}(t)= e^{kt}\psi(t)n_{i}(0)+\int\limits_{0}^{t} e^{k\left(t-t^{'}\right)}\psi(t-t^{'}) J_{i}^{+}(t^{'}) \; dt^{'}
\end{equation}
where the exponential factors  take into account the per capita addition or removal of particles 
as in  \cite{HLW2006}.

Now solving for the gain flux, $J^{+}_i(t)$, in  Eq.~(\ref{Model:III:Balance:lr}) we find 
\begin{equation}
\label{Model:III:Balance:rearr:lr}
 J_{i}^{+}(t)= J_{i}^{-}(t)+\frac{dn_i(t)}{dt} -k n_i(t)
\end{equation} 
and using Eq.~(\ref{Model:III:Loss:flux:lr}) gives
\begin{equation}
\label{Model:III:Loss:flux:sub:lr}
J_{i}^{-}(t)= e^{kt}\psi(t)n_{i}(0)+\int\limits_{0}^{t} e^{k\left(t-t^{'}\right)} \psi(t-t^{'}) \left[J_{i}^{-}(t^{'})+\frac{dn_i(t^{'})}{dt}-k n_i(t^{'})\right] \; dt^{'}.
\end{equation} 

Now using Laplace transform theory we find 
\begin{equation}
\label{Model:III:Loss:flux:lapl:lr}
\widehat{J}_{i}^{-}(s)= \widehat{\psi}(s-k)n_{i}(0)+ \widehat{\psi}(s-k) \left[\widehat{J}_{i}^{-}(s)+s\widehat{n}_i(s)-n_{i}(0)-k\widehat{n}_i(s) \right] 
\end{equation} 
and upon solving for the flux, we find a similar expression to Eq.~(\ref{Model:III:Loss:flux:lapl:2}),
\begin{equation}
\label{Model:III:Loss:flux:lapl:2:lr}
\widehat{J}_{i}^{-}(s)= \frac{\widehat{\psi}(s-k)}{\widehat{\Phi}(s-k)}\widehat{n}_i(s).
\end{equation} 
This Laplace trasnform can now be inverted
to find the loss flux given by the modified fractional derivative \cite{HLW2006,Sokolov2006} of the concentration at $i$,
\begin{equation}
\label{Model:III:Loss:flux:invert:lr}
{J}_{i}^{-}(t)= e^{kt}\frac{A_{\gamma}}{\tau^\gamma}\RLderivO{1-\gamma} \left(e^{-kt} n_{i}(t)\right),
\end{equation} 
where $A_\gamma$ is a constant given by $1$ or $1/\GammaFn{1-\gamma}$ if we use the waiting-time density in Eq.~(\ref{ModelII:Mittag:psi}) or Eq.~(\ref{Pareto:psi}) respectively.

Now using (\ref{Model:III:Gain:flux})  and (\ref{Model:III:Loss:flux:invert:lr}) in (\ref{Model:III:Balance:lr}) we find
\begin{multline}
\label{Model:III:Balance:2:lr}
\frac{dn_i(t)}{dt} = p_{r}(x_{i-1},t)
e^{kt}\frac{A_{\gamma}}{\tau^\gamma}\RLderivO{1-\gamma}\left(e^{-kt} n_{i-1}(t)\right)
+ p_{l}(x_{i+1},t)
e^{kt}\frac{A_{\gamma}}{\tau^\gamma}\RLderivO{1-\gamma}\left(e^{-kt} n_{i+1}(t)\right)\\-
e^{kt}\frac{A_{\gamma}}{\tau^\gamma}\RLderivO{1-\gamma}\left(e^{-kt} n_{i}(t)\right)
+k n_{i}(t).
\end{multline} 
The continuum limit following from setting
$x_{i}=x$,  $x_{i\pm 1}=x\pm \Delta x$, and Taylor series expansions about $x$, is given by
\begin{multline}
\label{ModelIII:eqn:final:lr}
\frac{\partial n}{\partial t} = D_{\gamma} e^{kt} \FracderivN{\gamma}  \left(e^{-kt} \frac{\partial^{2}n(x,t)}{\partial x^2}\right)
-\chi_{\gamma}\frac{\partial}{\partial x}\left(\frac{\partial c(x,t)}{\partial x}e^{kt} \FracderivN{\gamma} \left(e^{-kt} n(x,t)\right)\right)+kn(x,t).
\end{multline} 

\subsection{Model IV}
The master equation  for Model IV, Eq.(\ref{ModelIV:1}),  modified to include linear reaction dynamics is given by
 \begin{multline}
\label{ModelIV:1:lr}
n_{i}(t) =e^{kt} n_{i}(0) \Phi(t) + p_r(x_{i-1},t)\int\limits_0^{t}
 n_{i-1}(t^{'})  e^{k\left(t-t^{'}\right)}\psi(t-t^{'})\: dt^{'}
\\
+
p_l(x_{i+1},t)\int\limits_0^{t} n_{i+1}(t^{'}) e^{k\left(t-t^{'}\right)}\psi(t-t^{'}) \:dt^{'}.
\end{multline} 
Following similar steps used to simplify Eq.~(\ref{ModelIV:1}) and using Taylor series expansions we find
\begin{multline}
\label{ModelIV:4:lr}
\int\limits_0^{t} e^{k\left(t-t^{'}\right)}\Phi(t-t^{'}) \frac{\partial n(x,t^{'})}{\partial t} \:dt^{'}  \simeq  \frac{\Delta x^2}{2}
\int\limits_0^{t} \frac{\partial^2 n(x,t^{'})}{\partial x^2}   e^{k\left(t-t^{'}\right)}\psi(t-t^{'}) \:dt^{'}\\
+k\int\limits_0^{t} e^{k\left(t-t^{'}\right)}\Phi(t-t^{'}) n(x,t^{'}) \:dt^{'} 
-\Delta{x}^2\beta\frac{\partial}{\partial x}\left(\frac{\partial c(x,t)}{\partial x}\int\limits_0^{t} 
n(x,t^{'}) e^{k\left(t-t^{'}\right)}\psi(t-t^{'})\:dt^{'}\right)+O(\Delta x^4).
\end{multline}
Using Laplace transforms and the asymptotic expressions in Eqs.~(\ref{ModelII:alt:asym:Phi:lg}) and
(\ref{ModelII:alt:asym:psi:lg}) (evaluated for $s-k$ small) we arrive at Eq. (\ref{ModelII:eqn:final:lr}) for long times.

For short times we find 
\begin{equation}
\label{ModelII:4:alt:sh:lr}
\frac{\partial n(x,t)}{\partial t} =  D^{*}_\gamma e^{kt}
\FracIntN{1-\nu_\gamma}\left(e^{-kt} \frac{\partial^2 n(x,t)}{\partial x^2}\right) 
-\chi^{*}_\gamma e^{kt} \frac{\partial}{\partial x}\frac{\partial}{\partial t}\left(\frac{\partial c(x,t)}{\partial x} 
\FracIntN{-\nu_\gamma}\left(e^{-kt} n(x,t)\right)\right)+kn(x,t)
\end{equation} 
with $\nu_\gamma$  and the modified coefficients as defined previously for Model IV in section \ref{Models}.

\section{Numerical Solutions and Monte Carlo Simulations}
It is straightforward to obtain numerical solutions of the above model equations using difference approximations.
 In this section we describe  numerical solutions  for self-chemotactic variants of the model equations and we compare the solutions with Monte Carlo simulations.

Implementation details for the Monte Carlo simulations are described in the Appendix \ref{App:A}. In the results reported here simulations  were conducted on a one-dimensional lattice using the Pareto waiting-time density (\ref{Pareto:psi}) with the characteristic waiting-time $\tau=0.1$,  fractional exponent $\gamma=0.5$, and chemotactic sensitivity, $\beta$. The simulation results  are from an average of 200 runs with 10, 000 particles initially located at the origin. 

As our starting point for numerical solutions of the macroscopic models we consider the discrete space variants of Model II, III, and IV  given by Eq~(\ref{ModelII:8}), (\ref{Model:III:Balance:3}), and ~(\ref{ModelIV:1}) respectively.
A discrete space variant for Model I, analogous to the discrete space variant for Model II,  is given by
\begin{equation}
\label{Model:I:Meso}
\frac{dn_i(t)}{dt} =\frac{A_\gamma \gamma t^{\gamma-1}}{\tau^\gamma}\left[ p_{r}(x_{i-1},x_i,t) n_{i-1}(t)+
                     p_{l}(x_{i+1},x_i,t) n_{i+1}(t)-n_i(t)\right]
\end{equation} 

The discrete space equations for Models II and III were solved using an implicit time stepping method with the fractional derivatives approximated using the L1 scheme \cite{Oldham1974} as in \cite{Langlands2005}. For Model IV, the integrals in the discrete space representation were approximated by taking the unknown concentration, $n_{i\pm 1}(t)$, to be piecewise linear in time.

The numerical solutions of the discrete space equations and the Monte Carlo simulations were performed using the same space grid size and similar values for the parameters $\tau$,  $\gamma$, and $\beta$. We also chose the same initial condition (one at the origin and zero elsewhere). The constant $A_{\gamma}$ was chosen as in Eq.~(\ref{Agamma:Pd}) since we used the Pareto density, Eq.~(\ref{Pareto:psi}), in the simulations.

\begin{figure}
 \begin{center}
\includegraphics[angle=0,width=0.49\textwidth]{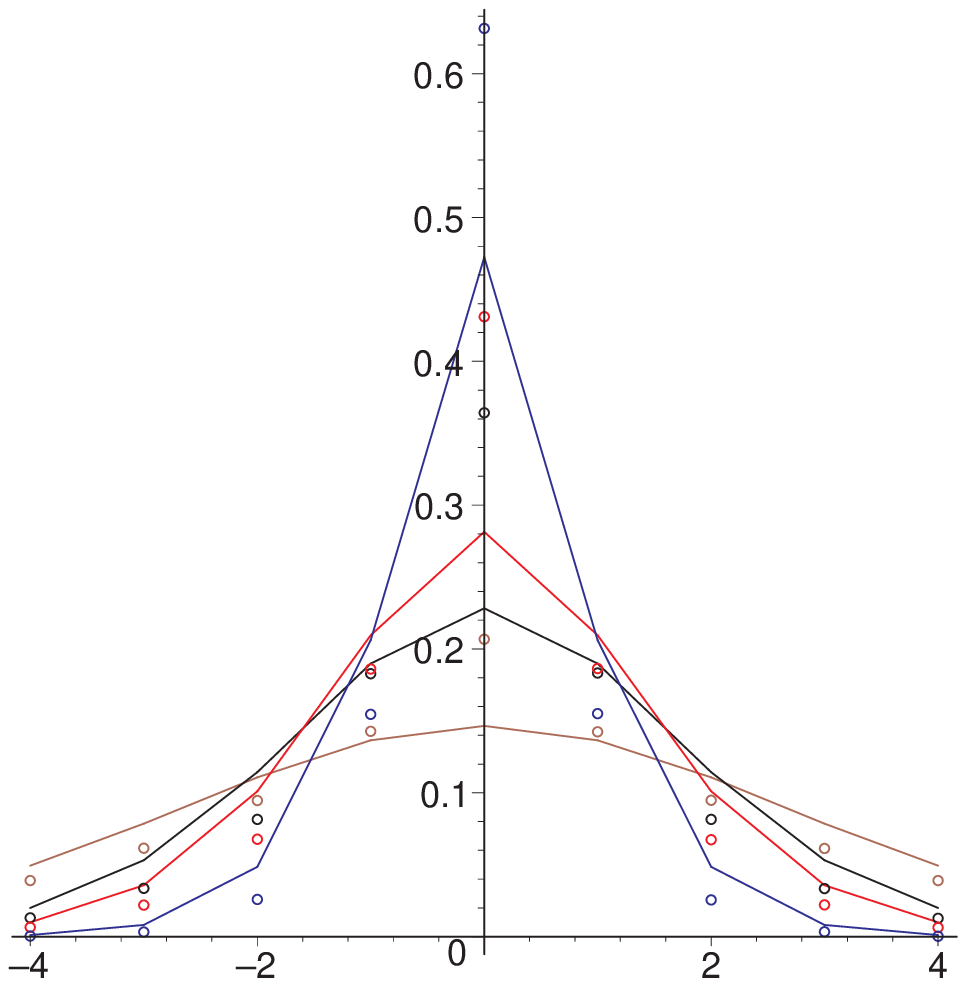}
 \includegraphics[angle=0,width=0.49\textwidth]{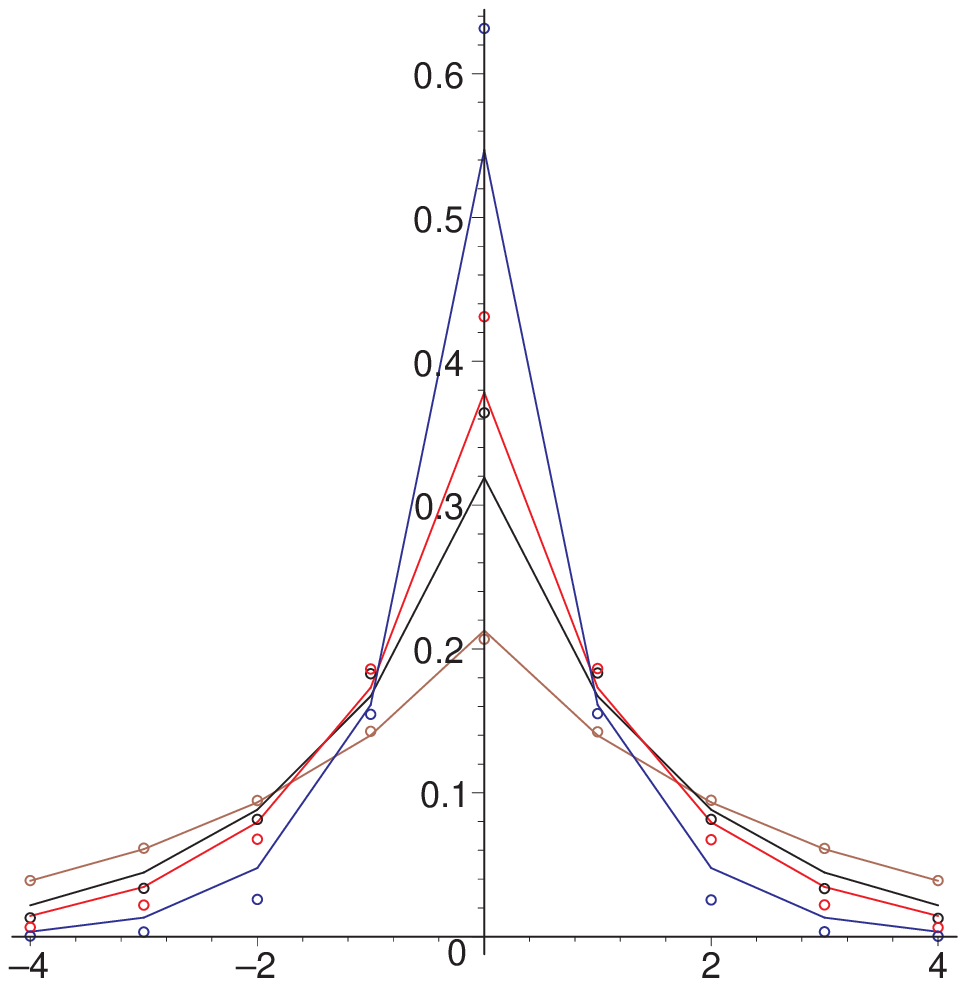}
 \includegraphics[angle=0,width=0.49\textwidth]{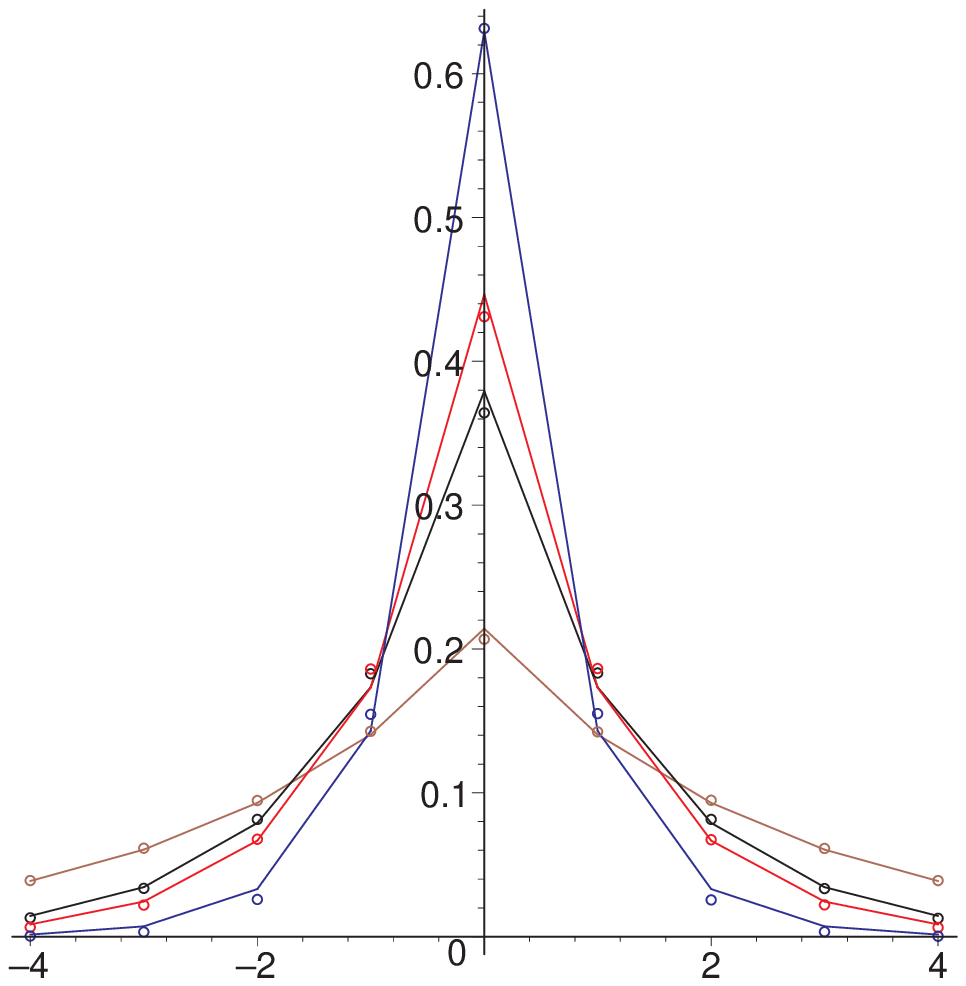}
 \includegraphics[angle=0,width=0.49\textwidth]{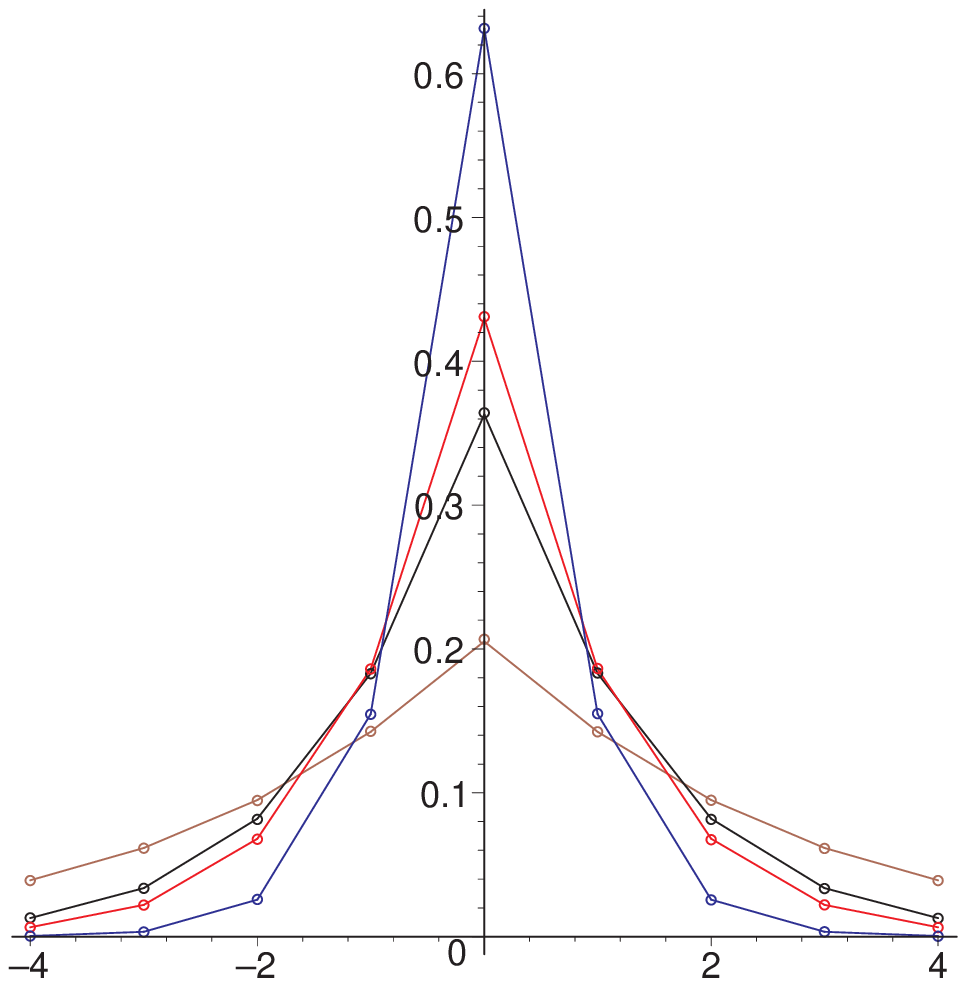}
\end{center}
\caption{\label{1D:beta:01:Prob:Comp} One-dimensional subdiffusive chemotaxis simulation results (circles) with $\beta=0.1$ showing the concentration profile at $t=0.4$ (blue), $t=2$ (red), $t=4$ (black), and $t=20$ (brown)  
compared with the numerical solution of the governing equations for models I~(top left), 
II~(top right), III~(bottom left) and IV~(bottom right) (solid lines). }
\end{figure}

\begin{figure}
 \begin{center}
\includegraphics[angle=0,width=0.49\textwidth]{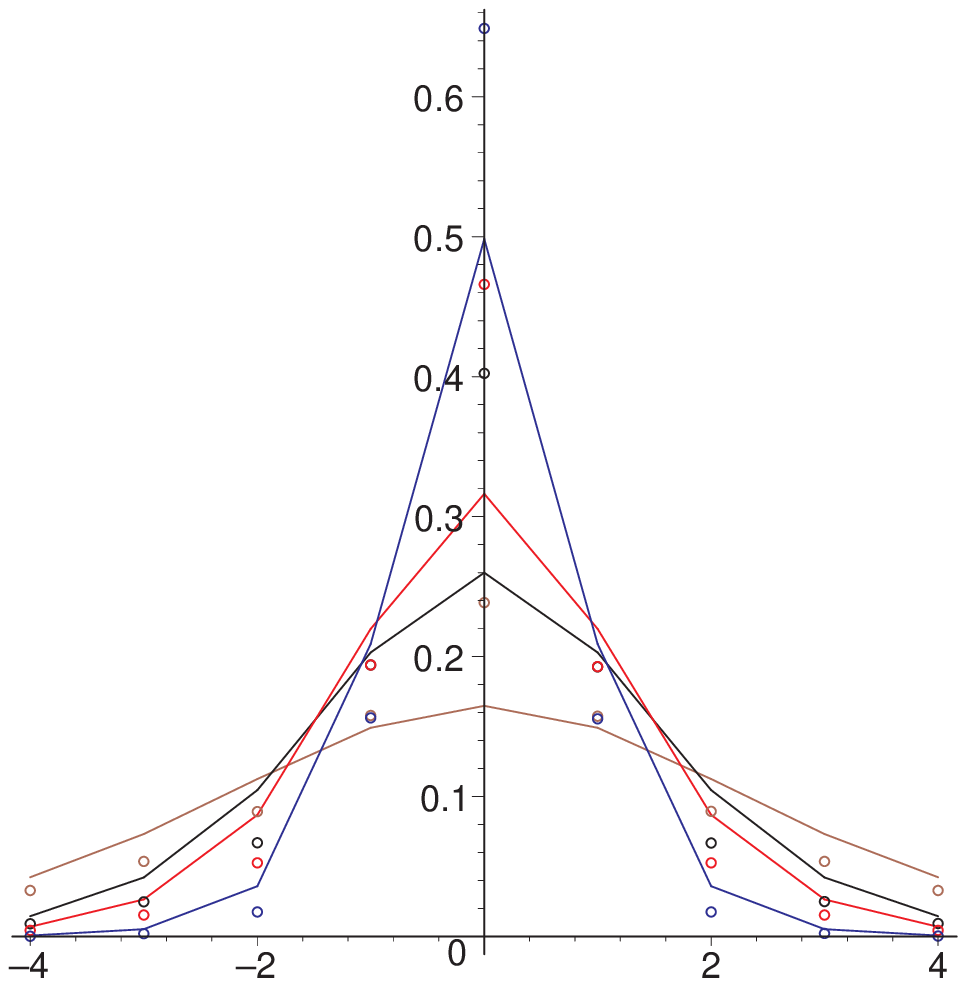}
 \includegraphics[angle=0,width=0.49\textwidth]{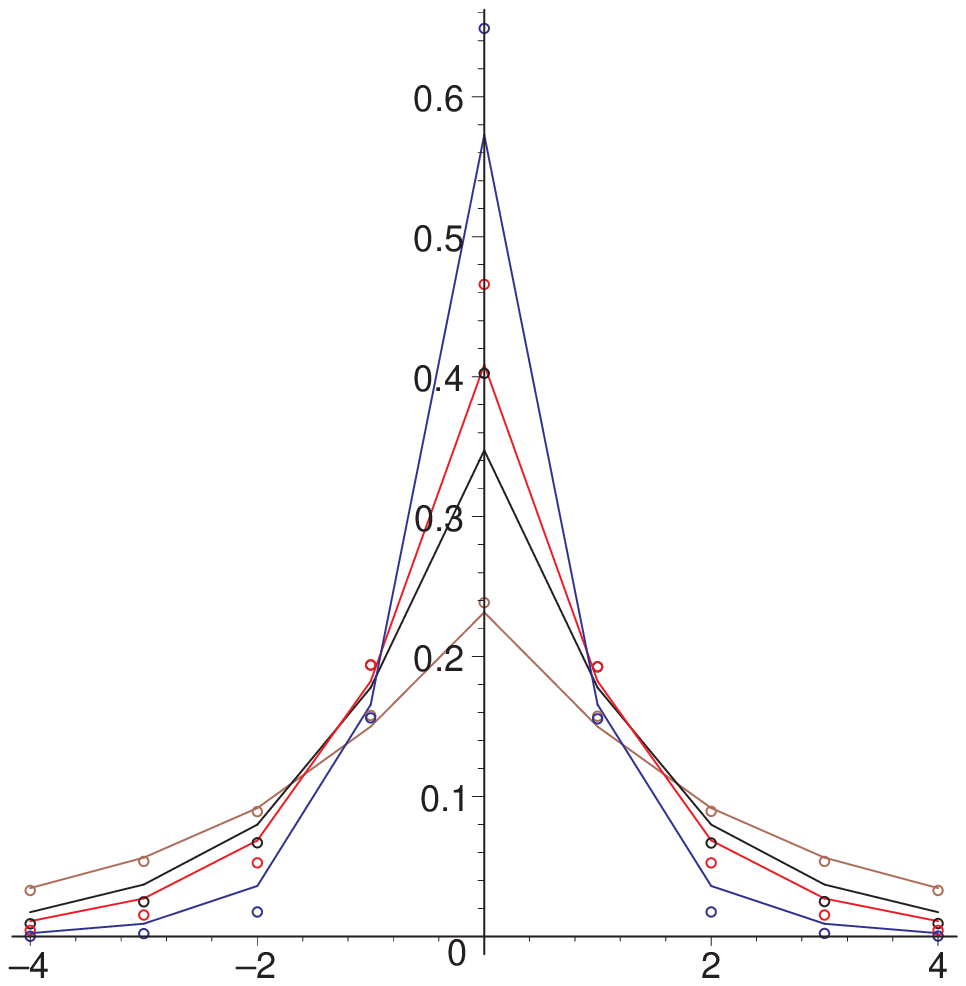}
 \includegraphics[angle=0,width=0.49\textwidth]{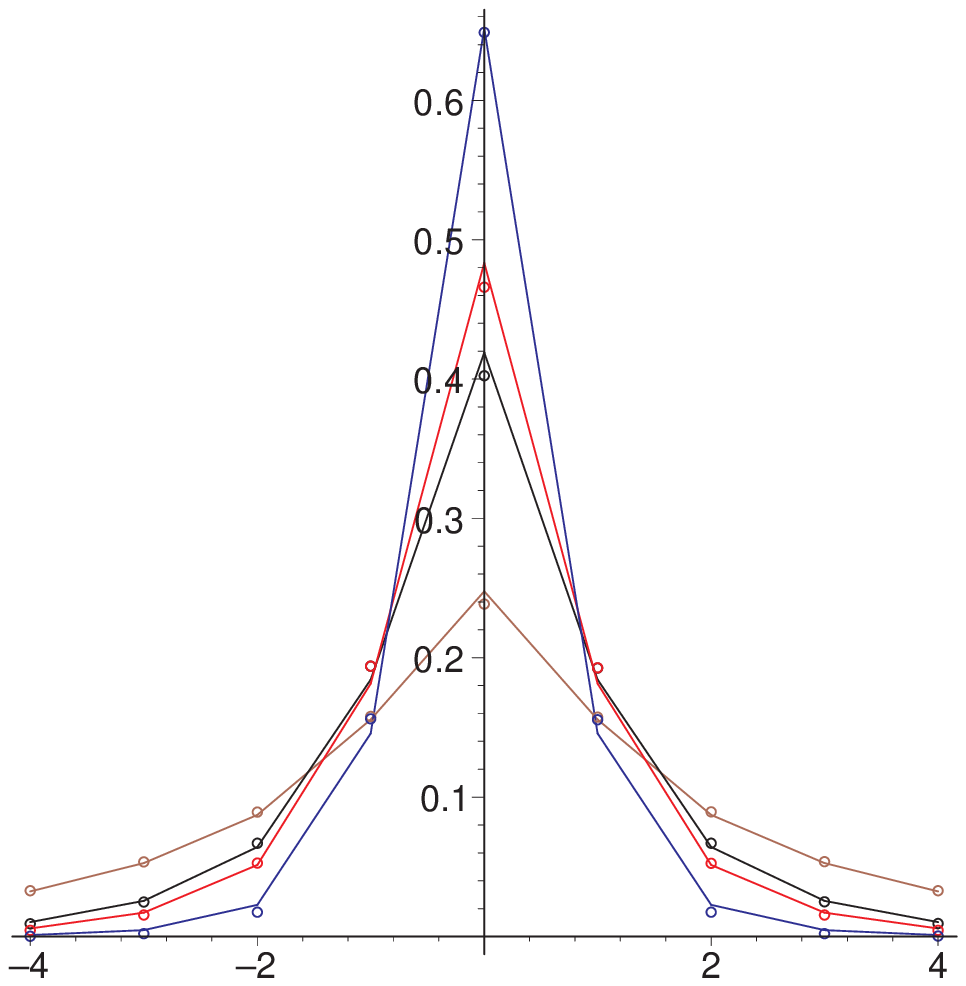}
 \includegraphics[angle=0,width=0.49\textwidth]{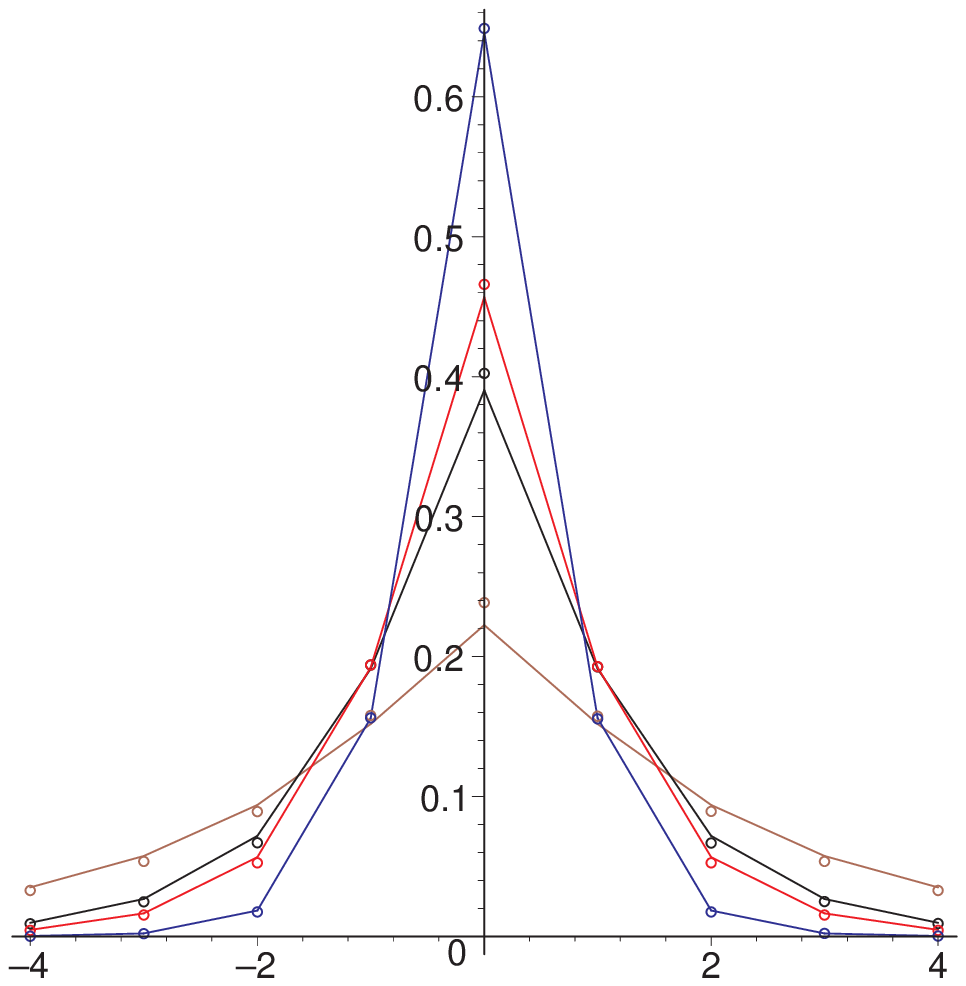}
\end{center}
\caption{\label{1D:beta:1:Prob:Comp} One-dimensional subdiffusive chemotaxis simulation results (circles) with $\beta=1.0$ showing the concentration profile at $t=0.4$ (blue), $t=2$ (red), $t=4$ (black), and $t=20$ (brown)  
compared with the numerical solution of the governing equations for models I~(top left), 
II~(top right), III~(bottom left) and IV~(bottom right) (solid lines). }
\end{figure}

\begin{figure}
 \begin{center}
\includegraphics[angle=0,width=0.49\textwidth]{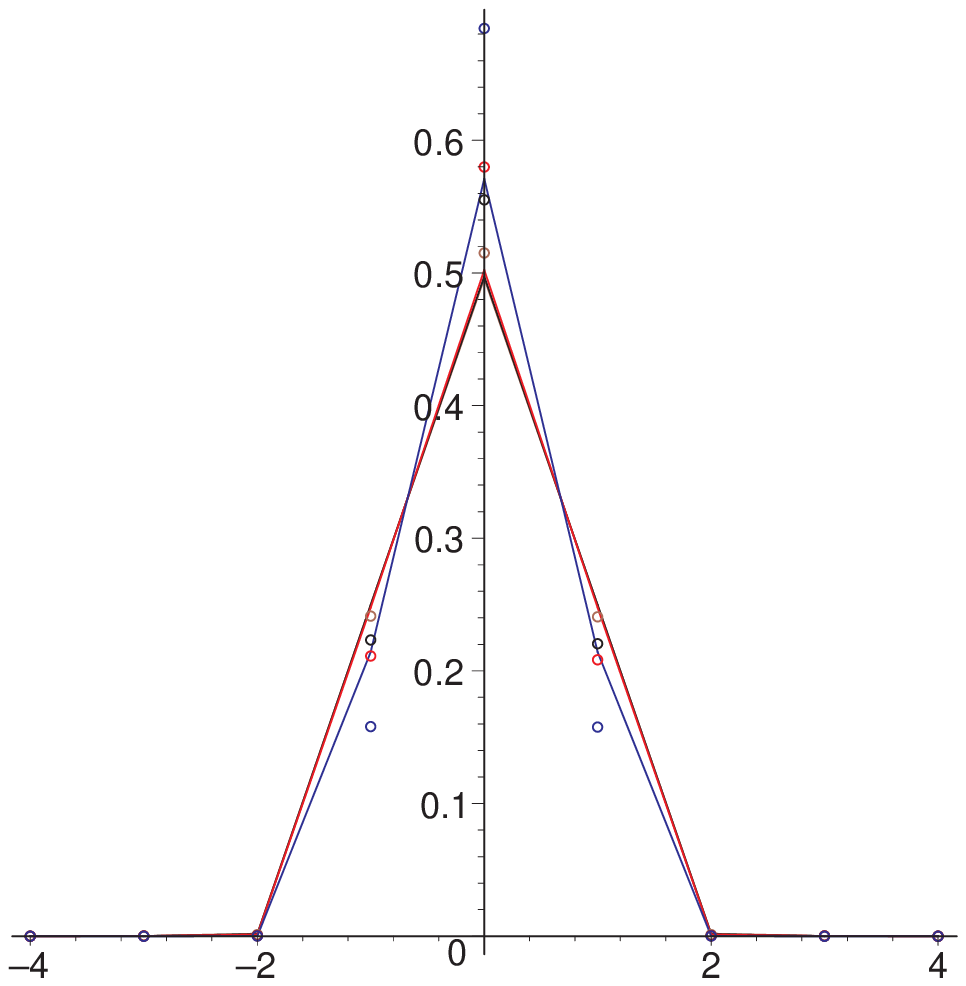}
 \includegraphics[angle=0,width=0.49\textwidth]{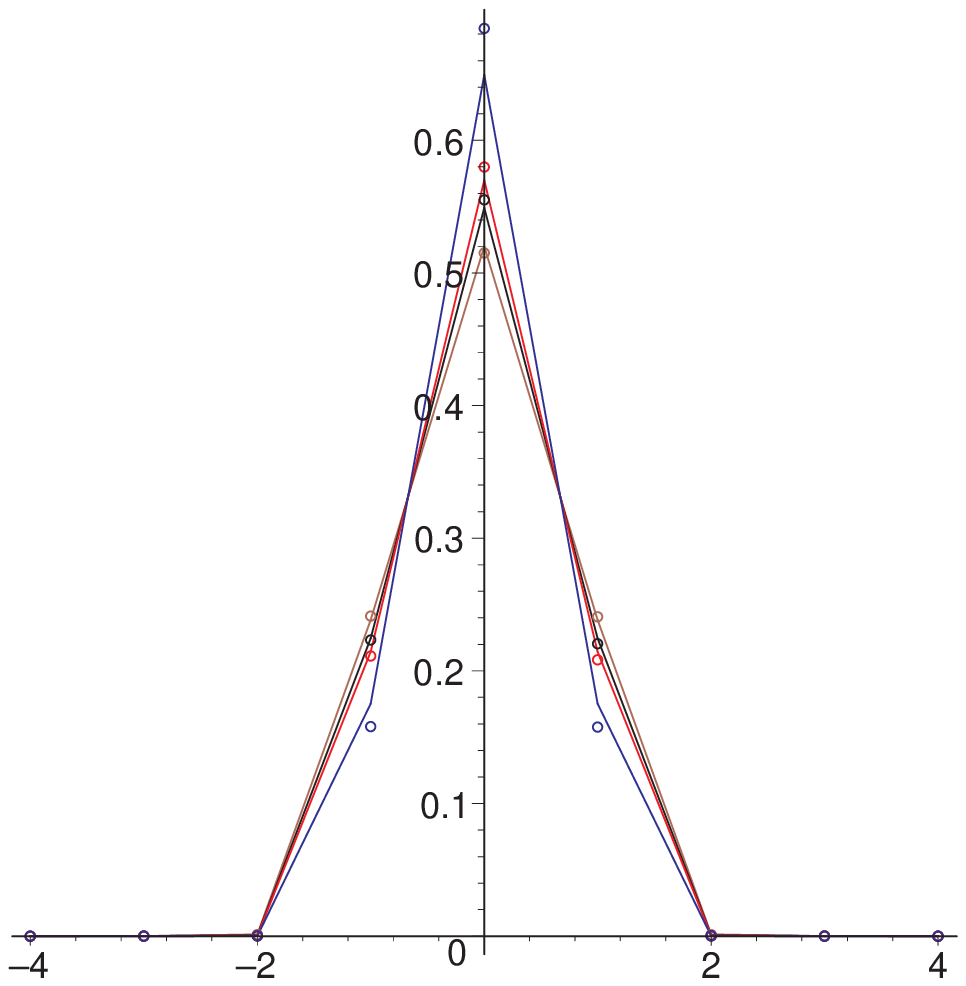}
 \includegraphics[angle=0,width=0.49\textwidth]{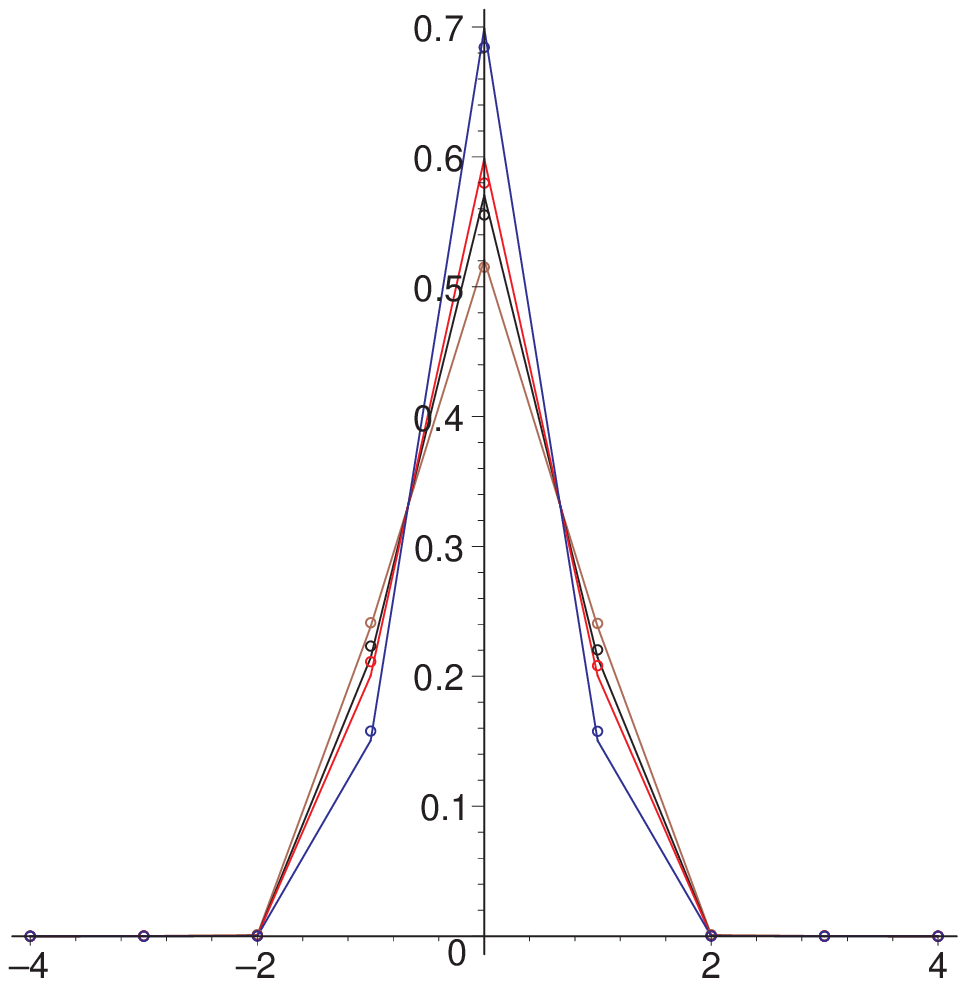}
 \includegraphics[angle=0,width=0.49\textwidth]{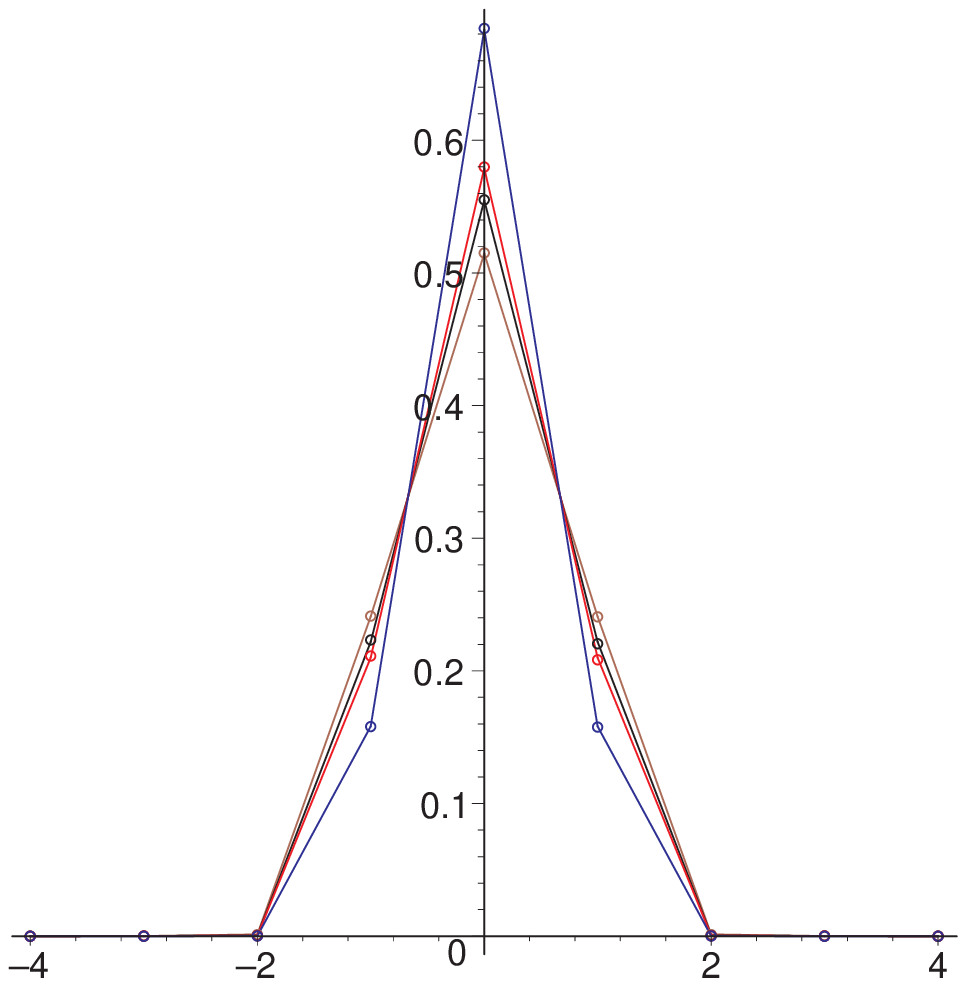}
\end{center}
\caption{\label{1D:beta:10:Prob:Comp} One-dimensional subdiffusive chemotaxis simulation results (circles) with $\beta=10.0$ showing the concentration profile at $t=0.4$ (blue), $t=2$ (red), $t=4$ (black), and $t=20$ (brown)  
compared with the numerical solution of the governing equations for models I~(top left), 
II~(top right), III~(bottom left) and IV~(bottom right) (solid lines). }
\end{figure}

In Figure~\ref{1D:beta:01:Prob:Comp} we compare the Monte Carlo simulation results 
with the numerical solution of the  discrete space equations for each model with the chemotactic sensitivity parameter, $\beta=0.1$. Further results are shown in
 Figures~\ref{1D:beta:1:Prob:Comp}  and \ref{1D:beta:10:Prob:Comp} for the sensitivity parameter values $\beta=1$ and $\beta=10$ respectively. 
 
 The numerical  solutions for Model III   Eq~(\ref{Model:III:Balance:3}) and Model IV (\ref{ModelIV:1}) are in close agreement with the Monte Carlo simulations at all times.
 The numerical solution for Model II does not fit with the Monte Carlo simulations well at short times but it provides a good fit at long times
 ($t=20$).
The numerical solution for Model I does not fit the Monte Carlo simulations well, especially for small values of $\beta$, where the predicted shape near the origin is smoother than that exhibited by the simulation data.

The closer agreement between the numerical results for
 Models III and IV  and the Monte Carlo simulations is due to the timing of the chemotactic forcing. In models III and IV, and in the Monte Carlo simulations the chemotactically influenced jumping probabilities are determined at the end of the waiting times, whereas in Model II they are determined at the start of the waiting times. This difference is less marked
 if the chemotactic concentration varies slowly in time.

Overall, the numerical solutions for Model III provide better agreement with the Monte Carlo simulations than the numerical solutions for Model IV. This better agreement can be seen at the intermediate value of $\beta=1$ in Figure~\ref{1D:beta:1:Prob:Comp}. This better agreement may be due to differences in numerical errors in approximating the discrete space equations for Model II and Model IV rather than due to differences between the equations themselves.


\section{Summary and Discussion}

The correct form of the fractional Fokker-Planck for particles undergoing anomalous subdiffusion in an external space and time varying force field has been an open problem. In the absence of a force field, subdiffusion can be modelled with a fractional temporal derivative operating on the spatial Laplacian. For subdiffusion in a purely space dependent force field the fractional temporal derivative can be put to the left of the standard terms on the right hand side of  the standard Fokker-Planck equation \cite{MBK1999,Barkai2000,MK2000,Sokolov2001} as in Eq.(\ref{fokkerplanck}). However the  consensus has been that  for subdiffusion in an external space-time-dependent force field the fractional temporal derivative should not operate on the force field \cite{HPGH2007,WMW2008,HPGH2009}, as in Eq.(\ref{dispersal}).
The modelling is further complicated if the force itself is affected directly or indirectly by the subdiffusing particles. This is the case in  fractional electro-diffusion and fractional chemotaxis diffusion, the case considered here.

In this article we have introduced and investigated four models for particles undergoing anomalous subdiffusion  in the presence of  chemotactic forcing. The first being based on an {\em adhoc} extension to the fractional Brownian motion equation (Model I), two models based on Continuous Time Random Walk (CTRW) master equations where concentration-dependent jump probabilities were evaluated before (Model II) or after (Model IV) the particle waiting, and
 a fourth model  derived from a  generalized master equation (Model III).  Concentration-dependent jump probabilities were used to incorporate the effect of chemotaxis in discrete space representations of the  models and in the Monte Carlo (MC) simulations. 

Evaluating the jump probabilities prior to waiting in 
the CTRW formulation (Model II) resulted in a macroscopic equation (valid in the long time limit) with the fractional derivative acting upon the chemotactic gradient. 
Conversely, using a generalized master equation approach with the probabilities evaluated after waiting but prior to jumping gave a macroscopic equation where the fractional derivative does not act upon the gradient (Model III).
The CTRW formulation  with the jump probabilities evaluated after waiting
(Model IV)  could only be reduced to a Fractional Fokker-Planck continuum equation in the asymptotic limit for long and short times.  For long-times Models II and IV coincide whilst for short-times we found Models III and IV coincide asymptotically if a Mittag-Leffler density is used.

We also introduced Monte Carlo methods for simulating anomalous sub-diffusion in a chemotactic force field. In the Monte Carlo simulations the chemotactically influenced jump lengths were computed at the end of the waiting times, similar to Models III and IV. This could explain
the excellent agreement we found between numerical solutions for Models III and IV and the Monte Carlo simulations. The numerical solutions for Model II also showed good agreement at long times. The numerical solutions based on the fractional Brownian motion equation, did not agree well with the Monte Carlo results.

The fractional chemotaxis diffusion models were further generalized to incorporate linear reaction dynamics. As in previous research, \cite{HLW2006,LHW2008,Sokolov2006} we found that the  incorporation of linear reactions required the replacement of the Riemann Liouville fractional derivative with a modified version, in addition to including the linear reaction term.

The fractional chemotaxis diffusion equations developed in this paper
provide a new class of models for biological transport  influenced by chemotactic forcing, macro-molecular crowding and traps.  We have recently generalized these models to include arbitrary space-and-time dependent forces \cite{HLS2010b}.

\begin{acknowledgements}
This work was supported by the Australian Research Council.
\end{acknowledgements}

\bibliographystyle{apsrev}
\bibliography{LH-PRE}

\appendix
\section{Monte Carlo Simulations}
\label{App:A}
In this section we briefly describe the Monte Carlo method used to simulate chemotaxis 
on a periodic one-dimensional lattice with long-tailed waiting-time density (subdiffusion).  
For each  simulation  run, the following steps are conducted 
\begin{enumerate}
\item Set the number of grid points, simulation time, and the initial number of particles.
\item Initialise the parameters for the waiting-time and jump-length probability density functions.
\item Set up the initial particle positions.
\item For each particle generate a random waiting-time, $\delta t$, (time of the first jump).
\item Initialise the output time  $t_{out}=\Delta{t}$.
\item \label{jump} Find the time of the next jump by finding the minimum of all jumping times, $t_{jump}$.
\item \label{test} If $t_{jump}>t_{out}$ then go to step (\ref{output}) otherwise go to step (\ref{surv}).
\item \label{output} Store the current particle positions. Add $\Delta{t}$ to $t_{out}$. If $t_{out}$ exceeds the simulation time then simulation ends otherwise go to step (\ref{test}). 
\item \label{surv} Generate a random jump-length, $\delta{x}$ (see below).
\item Generate a new waiting time and update this particle's jumping time and position
 ($t_{jump}=t_{jump}+\delta {t}$ $x_{new}=x_{old}+\delta{x}$ ).
\item Go to step (\ref{jump}) 
\end{enumerate}

\subsection{Generation of Waiting Times} 
The waiting times for each particle/jumper were generated by comparing a uniform random
number, $r\in (0,1)$, with the cumulative probability function of the waiting-time density.
We use the Pareto density (\ref{Pareto:psi}) as the density
which has used by Yuste, Acedo, and Lindeberg \cite{Yuste2004}. 
The generated waiting-time is given as
\begin{equation}
\delta t = \tau\left(\left(1-r\right)^{-\frac{1}{\gamma}}-1\right)
\end{equation} 
where $r\in (0,1)$ is a uniform random number.

\subsection{Generation of Jump Distances} 
The jump distance for each particle/jumper is generated by comparing an uniform random
number, $r\in (0,1)$, with the cumulative probability function of the jump-length density.

For the  simulations we use the jump-length probability density nearest neighbour jumps only:
\begin{equation}
\label{1d:jump}
\delta x = \left\lbrace\begin{array}{cl}
-\Delta{x}, & \: 0\le r < p_l,\\
 \Delta{x}, & \: p_l\le r < 1
\end{array}\right.
\end{equation} 
where $\Delta x$ is the grid spacing and $p_l=p_l(x_{i},t)$ is the probability of jumping to the left given previously in Eqs.~(\ref{ModelII:pleft}) and (\ref{ModelII:v}). 

To evaluate the probabilities of jumping to the left or right for
Eq.~(\ref{1d:jump}), requires the approximation of the chemoattractant concentration, $c(x_i,t)$, in Eq.~(\ref{ModelII:v}). This is estimated by the proportion of chemoattractant particles at the grid point, $x_i$, compared with the total number of particles in the system.

\end{document}